\documentclass[12pt]{amsart}

\usepackage{color,fancyhdr}
\usepackage[all]{xy}

\usepackage{amsmath,amsthm,mathrsfs,amsfonts,
amssymb,times,latexsym,mathabx}

\usepackage[usenames,dvipsnames,svgnames,table]{xcolor}

\DeclareMathAlphabet{\curly}{U}{rsfs}{m}{n}  

\newtheorem{theorem}{Theorem}[section]

\theoremstyle{definition}

\newtheorem{corollary}[theorem]{Corollary}

\theoremstyle{problem}
\newtheorem{problem}[theorem]{Problem}
\usepackage{xcolor}

\numberwithin{equation}{section}

\makeatletter
\renewcommand{\pmod}[1]{\allowbreak\mkern7mu({\operator@font mod}\,\,#1)}
\makeatother

\newcommand{\be}{\begin{equation}}
\newcommand{\ee}{\end{equation}}

\renewcommand{\le}{\leqslant}

\renewcommand{\ge}{\geqslant}

\begin{document}

\title[On disjoint sets]{On disjoint sets}

\author{Jin-Hui Fang*}
\address{Department of Mathematics, Nanjing University of Information Science $\&$ Technology, Nanjing 210044, PR China}
\email{fangjinhui1114@163.com}
\author{Csaba S\'{a}ndor}
\address{Department of Stochastics, Institute of Mathematics, Budapest University of Technology and Economics, M\H{u}egyetem rkp. 3., H-1111, Budapest, Hungary; Department of Computer Science and Information Theory, Budapest University of Technology and Economics, M\H{u}egyetem rkp. 3., H-1111 Budapest, Hungary; MTA-BME Lend\"{u}let Arithmetic Combinatorics Research Group, ELKH, M\H{u}egyetem rkp. 3., H-1111 Budapest, Hungary}
\email{csandor@math.bme.hu}
\thanks{* Corresponding author.}
\thanks{The first author is supported by the National Natural Science Foundation of China, Grant No. 12171246 and the Natural Science Foundation of Jiangsu Province, Grant No. BK20211282. The second author is supported by the NKFIH Grants No. K129335.}
\date{\today}%

\begin{abstract}
Two sets of nonnegative integers $A=\{a_1<a_2<\cdots\}$ and $B=\{b_1<b_2<\cdots\}$ are defined as \emph{disjoint}, if $\{A-A\}\bigcap\{B-B\}=\{0\}$, namely, the equation $a_i+b_t=a_j+b_k$ has only trivial solution. In 1984, Erd\H os and Freud [J. Number Theory 18 (1984), 99-109.] constructed disjoint sets $A,B$ with $A(x)>\varepsilon\sqrt{x}$ and $B(x)>\varepsilon\sqrt{x}$ for some $\varepsilon>0$, which answered a problem posed by Erd\H{o}s and Graham. In this paper, following Erd\H os and Freud's work, we explore further properties for disjoint sets. As a main result, we prove that, for disjoint sets $A$ and $B$, assume that $\{x_1<x_2<\cdots\}$ is a set of positive integers such that $\frac{A(x_n)B(x_n)}{x_n}\rightarrow 2$ as $x_n\to \infty$, then, (i) for any $0<c_1<c_2<1,$ $c_1x_n\le y\le c_2x_n$, we have $\frac{A(y)B(y)}{y}\rightarrow1$ as $n\rightarrow \infty$; (ii) for any $1<c_3<c_4<2,$ $c_3x_n\le y\le c_4x_n$, we have $A(y)B(y)=(2+o(1))x_n$ as $n\rightarrow \infty$.
\end{abstract}

\maketitle

\section{\bf Introduction}

As a classic topic, Sidon's problems attracts much interest in combinatorial number theory, where a set $A=\{a_k\}_{k=1}^{\infty}$ of positive integers is called \emph{Sidon} if the differences
$a_i-a_j(i\neq j)$ are all distinct. For a positive real number $x$, denote $A(x)$ by the number of elements $a\in A$ with $a\le x$. In \cite{Erdos2}, Erd\H os and Graham posed an interesting problem as follows: \vskip2mm

\noindent\textbf{Erd\H os-Graham Problem.} Let $A=\{a_1<a_2<\cdots\}$ and $B=\{b_1<b_2<\cdots\}$ be sequences of integers satisfying $A(x)>\varepsilon\sqrt{x}$ and $B(x)>\varepsilon\sqrt{x}$ for some $\varepsilon>0$. Is it true that $a_i-a_j=b_k-b_t$ has infinitely many solutions? \vskip2mm

In 1984, Erd\H os and Freud \cite{Erdos1} provided a negative answer by choosing $A$ as nonnegative integers consisting only even powers of two and $B$ as nonnegative integers consisting only odd powers of two, since every integer can be uniquely written as the sum of different powers of two and
\begin{eqnarray*}
\liminf_{x\rightarrow\infty}\frac{\min\{A(x)B(x)\}}{\sqrt{x}}=\frac{1}{\sqrt{2}}.
\end{eqnarray*}
In \cite{Erdos1}, Erd\H os and Freud also introduced the notations $SP$ and $IN$ as follows:
\begin{eqnarray*}
SP=\limsup_{x\rightarrow\infty}\frac{A(x)B(x)}{x},\hskip4mm IN=\liminf_{x\rightarrow\infty}\frac{\min\{A(x),B(x)\}}{\sqrt{x}}.\end{eqnarray*}
Two sets of nonnegative integers $A=\{a_1<a_2<\cdots\}$ and $B=\{b_1<b_2<\cdots\}$ are defined as \emph{disjoint}, if $\{A-A\}\bigcap\{B-B\}=\{0\}$, namely, the equation $a_i+b_t=a_j+b_k$ has only trivial solution. Erd\H os and Freud \cite{Erdos1} proved the following nice results for disjoint sets $A$ and $B$:\vskip2mm

\noindent\textbf{Theorem A.} (\cite[Theorem 1]{Erdos2})
The largest possible value of $SP$ is 2, moreover the following more precise estimations hold:\vskip2mm

\noindent 1.1. To any function $H(x)$ with $\limsup_{x\rightarrow\infty}H(x)=\infty$, we can construct sets $A,B$ such that
 \begin{eqnarray*}A(x)B(x)\ge 2x-H(x)\end{eqnarray*}
is valid for infinitely many (integer) values of $x$.\vskip2mm

\noindent 1.2. The previous result is best possible: for any $A$ and $B$, $A(x)B(x)-2x\rightarrow-\infty$ ($x\rightarrow\infty$).\vskip2mm

\noindent\textbf{Theorem B.} (\cite[Theorem 3.3]{Erdos2})
If $SP=2$, then $IN=0$.\vskip2mm

In this paper, we further consider the properties of disjoint sets and prove that:\vskip2mm

\begin{theorem}\label{thm1} For disjoint sets $A$ and $B$, assume that $\{x_1<x_2<\cdots\}$ is a set of positive integers such that \begin{eqnarray*}\lim_{n\rightarrow \infty}
\frac{A(x_n)B(x_n)}{x_n}=2, \hskip3mm \mbox{namely,}\hskip3mm SP=2.\end{eqnarray*} \vskip2mm

\noindent (i) For any $0<c_1<c_2<1,$ $c_1x_n\le y\le c_2x_n$, we have
\begin{eqnarray*}\frac{A(y)B(y)}{y}\rightarrow1 \hskip3mm\mbox{as}\hskip3mm n\rightarrow \infty.
\end{eqnarray*}\vskip2mm

\noindent (ii) For any $1<c_3<c_4<2,$ $c_3x_n\le y\le c_4x_n$, we have
\begin{eqnarray*}A(y)B(y)=(2+o(1))x_n \hskip3mm\mbox{as}\hskip3mm n\rightarrow \infty.
\end{eqnarray*}\vskip2mm

\noindent (iii) \begin{eqnarray*}\limsup_{n\rightarrow \infty}\frac{A(2x_n)B(2x_n)}{2x_n}\le \frac{3}{2}.\end{eqnarray*}
Furthermore, the constant $\frac{3}{2}$ is best possible.
\end{theorem}\vskip2mm

We characterize the set $\{x_1<x_2<\cdots\}$ such that $\frac{A(x_n)B(x_n)}{x_n}\rightarrow 2$ for disjoint sets $A$ and $B$.\vskip2mm

\begin{theorem}\label{thm2} (i) For disjoint sets $A$ and $B$, if $\{x_1<x_2<\cdots\}$ is a set of positive integers with $x_{i+1}\ge 2x_i$ for all large $i$ such that \begin{eqnarray}\label{08191}\lim_{n\rightarrow \infty}
\frac{A(x_n)B(x_n)}{x_n}=2,\end{eqnarray} then \begin{eqnarray}\label{08172}\frac{x_{n+1}}{x_n}\rightarrow\infty\hskip4mm \mbox{as}\hskip4mm n\rightarrow\infty.\end{eqnarray}\vskip1mm

\noindent (ii) Let $\{x_1<x_2<\cdots\}$ be a set of positive integers satisfying \eqref{08172}. Then there exist disjoint sets $A$ and $B$ such that \begin{eqnarray*}\label{03195}\lim_{n\rightarrow \infty}\frac{A(x_n)B(x_n)}{x_n}=2.\end{eqnarray*}
\end{theorem}

Furthermore, we generalize the above Theorem B, that is:\vskip2mm

\begin{theorem}\label{thm3} For disjoint sets $A$ and $B$, if $IN\ge0$, then \begin{eqnarray}\label{thm201} \frac{IN^6}{IN^2+8}\le 64(2-SP).\end{eqnarray}
\end{theorem}\vskip2mm

For disjoint sets $A$ and $B$, we have $\binom{A(x)}{2}+\binom{B(x)}{2}\le x$. Hence $IN\le 1$.
The following corollary is obtained.\vskip2mm

\begin{corollary}  For disjoint sets $A$ and $B$, we have $IN\le 2\sqrt[3]{3}\sqrt[6]{2-SP}$.
\end{corollary}
For a given positive integer $k\ge2$, write
\begin{eqnarray*}
&&A=\{\epsilon_0+\epsilon_2k^2+\cdots+\epsilon_{2s-2}k^{2s-2}, \epsilon_{2i}=0,1,\cdots,k-1\},\nonumber\\
&&B=\{\epsilon_1k+\epsilon_3k^3+\cdots+\epsilon_{2s-1}k^{2s-1}, \epsilon_{2i-1}=0,1,\cdots,k-1\},
\end{eqnarray*} Obviously, $A$ and $B$ are disjoint sets. Erd\H os and Freud obtained in \cite{Erdos1} that
\begin{eqnarray*}SP=\frac{2(k+1)}{k+2} \hskip3mm\mbox{and}\hskip3mm IN=\frac{1}{\sqrt{k}}.\end{eqnarray*} Based on the above result, we naturally posed a problem:\vskip2mm
\begin{problem}  For disjoint sets $A$ and $B$, is there some positive constant $C$ such that
$$IN\le C\sqrt{2-SP}?$$
\end{problem}

\vskip3mm

\section{\bf Proof of Main Results}\vskip2mm

During the proof of Theorems \ref{thm1}-\ref{thm2}, we will employ the following fact many times.\vskip3mm

\noindent\textbf{Proposition 1.} Let $A$ and $B$ be disjoint sets with $SP=2$, namely, \begin{eqnarray*}\lim_{n\rightarrow \infty}\frac{A(x_n)B(x_n)}{x_n}=2.\end{eqnarray*} Then, with the exception of at most $o(x_n)$ numbers, all numbers in $[0,2x_n]$ can be written as $a_i+b_j$ with $a_i\le x_n, b_j\le x_n$.\vskip2mm

\noindent{\bf Proof of Theorem \ref{thm1}.} Before the proof of Theorem \ref{thm1}, we first prove the following proposition.\vskip2mm

\noindent\textbf{Proposition 2.} There exist functions $y_n<z_n<x_n$, $y_n \to \infty$ such that $y_n=o(x_n)$, $z_n=x_n-y_n$ and one of the following conclusion holds:\vskip2mm

\noindent(i) \begin{eqnarray*}&&A(y_n)=(0.5+o(1))A(x_n),\hskip2mm A(z_n)-A(y_n)=o(x_n),\hskip2mm A(x_n)-A(z_n)=(0.5+o(1))A(x_n),\end{eqnarray*} and \begin{eqnarray*}B(y)=(\frac{y}{x_n}+o(1))B(x_n)\hskip3mm \mbox{for}\hskip3mm y_n<y<z_n;\end{eqnarray*}\vskip2mm

\noindent(ii) \begin{eqnarray*}&&B(y_n)=(0.5+o(1))B(x_n),\hskip2mm B(z_n)-B(y_n)=o(x_n),\hskip2mm B(x_n)-B(z_n)=(0.5+o(1))B(x_n)\end{eqnarray*} and \begin{eqnarray*}A(y)=(\frac{y}{x_n}+o(1))A(x_n)\hskip3mm \mbox{for}\hskip3mm y_n<y<z_n;\end{eqnarray*}\vskip2mm

\begin{proof} For any positive number $c'$ and any $y_n$ with $c'x_n\le y_n\le \frac{1}{2}x_n$, we can only use $a_i\le y_n$ and $b_i\le y_n$ for the numbers in $[0,y_n]$. Similarly, we can only use $z_n\le a_i\le x_n$ and $z_n\le b_i\le x_n$ for the numbers in $[u_n,2x_n]$. It infers from Proposition 1 that
\begin{eqnarray}\label{101}
A(y_n)B(y_n)>y_n-o(x_n),\hskip3mm (A(x_n)-A(z_n)(B(x_n)-B(z_n))>y_n-o(x_n).\end{eqnarray}
On the other hand, by considering the difference, we have $$A(y_n)B(y_n)+(A(x_n)-A(z_n))(B(x_n)-B(z_n))\le 2y_n.$$ It follows from \eqref{101} that
\begin{eqnarray}\label{102}
A(y_n)B(y_n)=y_n+o(x_n)\end{eqnarray}
and
\begin{eqnarray}\label{103}
(A(x_n)-A(z_n))(B(x_n)-B(z_n))=y_n+o(x_n)\end{eqnarray}

Let $\delta_1>0$ and let $\delta_1 x_n\le y_n\le\frac{1}{3}x_n$. Then the sets $A\cap [0,y_n]$, $A\cap ]y_n,2y_n]$ and $A\cap]z_n,x_n]$ are distinct. By considering the differences we have $$A(y_n)B(y_n)+(A(2y_n)-A(y_n))(B(2y_n)-B(y_n))+(A(x_n)-A(z_n))(B(x_n)-B(z_n))\le 2y_n.$$ It follows from (\ref{102}) and (\ref{103}) that
\begin{eqnarray}\label{104}
(A(2y_n)-A(y_n))(B(2y_n)-B(y_n))=o(x_n).\end{eqnarray} Since $y_n\ge \delta_1 x_n$, it follows from (\ref{102}) that $A(y_n)B(y_n)\ge \delta_1 x_n+o(x_n)$. It infers from (\ref{104}) that
\begin{eqnarray}\label{107} A(2y_n)-A(y_n)=o(A(y_n)) \hskip3mm\mbox{or}\hskip3mm B(2y_n)-B(y_n)=o(B(y_n)). \end{eqnarray}
Similarly, by considering the sets $A\cap [0,y_n]$, $A\cap ]z_n,x_n-2y_n]$ and $A\cap]z_n,x_n]$ we get \begin{eqnarray}\label{110} &&\hskip6mmA(z_n)-A(x_n-2y_n)=o(A(x_n)-A(z_n)) \\ &&\nonumber\mbox{or}\hskip3mm B(z_n)-B(x_n-2y_n)=o(B(x_n)-B(z_n)). \end{eqnarray}
It follows from (\ref{102}), (\ref{103}) and the fact $A(x_n)B(x_n)=2x_n+o(x_n)$ that \begin{eqnarray*}A(x_n)B(y_n)+A(y_n)B(x_n)=2x_n+o(x_n).\end{eqnarray*}
Then, \begin{eqnarray*}\frac{A(x_n)B(y_n)+A(y_n)B(x_n)}{A(x_n)B(x_n)}
=\frac{B(y_n)}{B(x_n)}+\frac{A(y_n)}{A(x_n)}=1+o(1).\end{eqnarray*}

Take $y_n=\frac{x_n}{2}$. Let $\alpha_n=\frac{A(\frac{x_n}{2})}{A(x_n)}$ and $\beta_n=\frac{B(\frac{x_n}{2})}{B(x_n)}$. Clearly,  $\alpha_n \beta_n =0.25+o(1)$ and $\alpha_n +\beta_n =1+o(1)$, that is $\alpha_n=0.5+o(1)$ and $\beta_n =0.5+o(1)$. Hence \begin{eqnarray}\label{106} A(\frac{x_n}{2})=(0.5+o(1))A(x_n) \hskip3mm\mbox{and}\hskip3mm B(\frac{x_n}{2})=(0.5+o(1))B(x_n). \end{eqnarray}
\vskip2mm

Let us suppose that \begin{eqnarray}\label{108} A(\frac{x_n}{2})-A(\frac{x_n}{4})=o(A(\frac{x_n}{2})). \end{eqnarray} We will prove by induction on $k$ that \begin{eqnarray}\label{08211} A(\frac{x_n}{2})-A(\frac{x_n}{2^k+1})=o(A(\frac{x_n}{2})). \end{eqnarray} for every positive integer $k$.\vskip2mm

It follows from (\ref{108}) that $A(\frac{x_n}{2})-A(\frac{x_n}{3})=o(A(\frac{x_n}{2}))$, thus \eqref{08211} holds for $k=1$. Suppose that \begin{eqnarray*}\label{109} A(\frac{x_n}{2})-A(\frac{x_n}{2^k+1})=o(A(\frac{x_n}{2})).\end{eqnarray*} Then \begin{eqnarray*}A(\frac{x_n}{2^k+1})=(1+o(1))A(\frac{x_n}{2}).\end{eqnarray*}
Since
\begin{eqnarray*}[0,\frac{x_n}{2^k+1}]\subseteq [0,2\frac{x_n}{2^{k+1}+1}]
\subseteq [0,\frac{x_n}{2}],\end{eqnarray*}
we have \begin{eqnarray}\label{08221}A(2\frac{x_n}{2^{k+1}+1})=(1+o(1))A(\frac{x_n}{2}).\end{eqnarray}
Thus, \begin{eqnarray}\label{08215}
A(2\frac{x_n}{2^{k+1}+1})=(1+o(1))A(\frac{x_n}{2^k+1}).
\end{eqnarray}
Take $y_n=\frac{x_n}{2^{k+1}+1}$ in (\ref{107}), then
\begin{eqnarray}\label{08212}&&\hskip6mmA(2\frac{x_n}{2^{k+1}+1})-A(\frac{x_n}{2^{k+1}+1})
=o(A(\frac{x_n}{2^{k+1}+1}))\\
&&\nonumber\mbox{or} \hskip2mm B(2\frac{x_n}{2^{k+1}+1})-B(\frac{x_n}{2^{k+1}+1})=o(B(\frac{x_n}{2^{k+1}+1})).\end{eqnarray}
Assume that $B(2\frac{x_n}{2^{k+1}+1})-B(\frac{x_n}{2^{k+1}+1})=o(B(\frac{x_n}{2^{k+1}+1}))$. Then \begin{eqnarray*}B(2\frac{x_n}{2^{k+1}+1})-B(\frac{x_n}{2^k+1})=o(B(\frac{x_n}{2^k+1})) \hskip3mm \mbox{namely,}\hskip3mm B(2\frac{x_n}{2^{k+1}+1})=(1+o(1))B(\frac{x_n}{2^k+1}).\end{eqnarray*} It follows from (\ref{102}) and \eqref{08215} that $$(1+o(1))2\frac{x_n}{2^{k+1}+1}=A(2\frac{x_n}{2^{k+1}+1})B(2\frac{x_n}{2^{k+1}+1})=(1+o(1))A(\frac{x_n}{2^k+1})B(\frac{x_n}{2^k+1})=(1+o(1))\frac{x_n}{2^k+1},$$
a contradiction. Then we could deduce from \eqref{08212} that \begin{eqnarray*}A(2\frac{x_n}{2^{k+1}+1})-A(\frac{x_n}{2^{k+1}+1})=o(A(\frac{x_n}{2^{k+1}+1})).
\end{eqnarray*} It follows from \eqref{08221} that $$A(\frac{x_n}{2})-A(\frac{x_n}{2^{k+1}+1})=o(A(\frac{x_n}{2})),$$ which completes the induction.
Thus, \eqref{08211} follows.
\vskip2mm

If \eqref{08211} holds, then $A(\delta x_n)=(1+o(1))A(\frac{x_n}{2})$ for every $0<\delta\le \frac{1}{2}$. We will show that: if $A(\delta_2 x_n)=(1+o(1))A(\frac{x_n}{2})$ for every $0<\delta_2\le \frac{1}{2}$ and $A(Cx_n)=(1+o(1))A(\frac{x_n}{2})$ for some $C$ with $0<C<1$, then $B(cx_n)=(c+o(1))B(x_n)$ for every $0<c\le C$. Let $0<\delta <c$. Clearly,
by \eqref{101} we know that $$
A(cx_n)B(cx_n)\ge cx_n-o(x_n).
$$
On the other hand,
\begin{eqnarray*}
A(cx_n)B(cx_n)&=&(A(cx_n)-A(\delta x_n))B(cx_n)+A(\delta x_n)B(cx_n)\\
&=&o(A(\frac{x_n}{2}))B(cx_n)+|\{ m: m=a+b, a\le \delta x_n, a\in A, b\le cx_n, b\in B\}|\\
&\le& (c+\delta +o(1))x_n.\end{eqnarray*}
Since $\delta$ can be arbitrary small, we have $A(cx_n)B(cx_n)=(c+o(1))x_n$. Then
we could deduce from \eqref{102} (choosing $y=\frac{x_n}{2}$) and \eqref{106} that
$$A(cx_n)B(cx_n)=(c+o(1))x_n=(2c+o(1))A(\frac{x_n}{2})B(\frac{x_n}{2})=(c+o(1))A(cx_n)B(x_n).
$$
It follows that $B(cx_n)=(c+o(1))B(x_n)$ for any $c$ with $0<c\le C$. Since $A(\frac{x_n}{2})=(1+o(1))A(\frac{x_n}{2})$, we have $B(cx_n)=(c+o(1))B(x_n)$ for any $c$ with $0<c\le 0.5$.\vskip2mm

Next, we will prove that \begin{eqnarray}\label{111} A(x_n-\frac{x_n}{2^k+1})=(1+o(1))A(\frac{x_n}{2}). \end{eqnarray} We argue by induction on $k$. Taking $y=\frac{x_n}{3}$ in (\ref{110}) we get
\begin{eqnarray}\label{08213}
&&\hskip6mm A(x_n-\frac{x_n}{3})-A(x_n-2\frac{x_n}{3})=o(A(x_n)-A(\frac{2x_n}{3}))\\
&&\nonumber\mbox{or} \hskip2mm B(x_n-\frac{x_n}{3})-B(x_n-2\frac{x_n}{3})=o(B(x_n)-B(\frac{2x_n}{3})).\end{eqnarray}
If $B(\frac{2x_n}{3})-B(\frac{x_n}{3})=o(B(x_n)-B(\frac{2x_n}{3}))$, then \begin{eqnarray*}B(\frac{x_n}{2})-B(\frac{x_n}{3})=o(B(x_n))=o(B(\frac{x_n}{2})).\end{eqnarray*} Hence $B(\frac{x_n}{3})=(1+o(1))B(\frac{x_n}{2})$. But $B(\frac{x_n}{3})=(\frac{1}{3}+o(1))B(x_n)$ and $B(\frac{x_n}{2})=(\frac{1}{2}+o(1))B(x_n)$, a contradiction. Then we could deduce from \eqref{08213} that
\begin{eqnarray*}A(\frac{2x_n}{3})-A(\frac{x_n}{3})
=o(A(x_n)-A(\frac{2x_n}{3}))=o(A(\frac{x_n}{2})).\end{eqnarray*} It follows that $A(\frac{2x_n}{3})=(1+o(1))A(\frac{x_n}{2})$, thus \eqref{111} holds for $k=1$. Suppose that \begin{eqnarray}\label{08216}A(x_n-\frac{x_n}{2^k+1})=(1+o(1))A(\frac{x_n}{2}).\end{eqnarray} Noting that
\begin{eqnarray*} \frac{x_n}{2}<x_n-2\frac{x_n}{2^{k+1}+1}<x_n-\frac{x_n}{2^k+1},\end{eqnarray*}
we know that $A(x_n-2\frac{x_n}{2^{k+1}+1})=(1+o(1))A(\frac{x_n}{2})$. Taking $y=\frac{x_n}{2^{k+1}+1}$ in (\ref{110}) we get \begin{eqnarray}\label{08214}
&&\hskip6mmA(x_n-\frac{x_n}{2^{k+1}+1})-A(x_n-2\frac{x_n}{2^{k+1}+1})
=o(A(x_n)-A(x_n-\frac{x_n}{2^{k+1}+1}))\\
&&\nonumber\mbox{or} \hskip2mm
B(x_n-\frac{x_n}{2^{k+1}+1})-B(x_n-2\frac{x_n}{2^{k+1}})=o(B(x_n)-B(\frac{x_n}{2^{k+1}+1})).
\end{eqnarray}
If $B(x_n-\frac{x_n}{2^{k+1}+1})-B(x_n-2\frac{x_n}{2^{k+1}+1}))
=o(B(x_n)-B(x_n-\frac{x_n}{2^{k+1}+1}))$, then \begin{eqnarray*}B(x_n-\frac{x_n}{2^k+1})-B(x_n-2\frac{x_n}{2^{k+1}+1})=o(B(x_n)).
\end{eqnarray*} Hence
\begin{eqnarray*}
&&(\frac{1}{2(2^k+1)(2^{k+1}+1)}+o(1))x_n\\
&\le&|\{ m:m=a+b, a\in A, b\in B,x_n-\frac{x_n}{2^k+1}+\frac{1}{2(2^k+1)(2^{k+1}+1)}x_n\le m
\le x_n-2\frac{x_n}{2^{k+1}+1}\} |\\
&\le&| \{ m:m=a+b, a\in A, a\le \frac{x_n}{2(2^k+1)(2^{k+1}+1)},\\ && x_n-\frac{x_n}{2^k+1}+\frac{x_n}{2(2^k+1)(2^{k+1}+1)}
\le m\le x_n-2\frac{x_n}{2^k+1}\} |\\ &&+|\{ m: m=a+b, a\in A, a> \frac{x_n}{2(2^k+1)(2^{k+1}+1)},\\
 &&x_n-\frac{x_n}{2^k+1}+\frac{1}{2(2^k+1)(2^{k+1}+1)}x_n\le m\le x_n-2\frac{x_n}{2^k+1}\}|\\
&\le& A(\frac{x_n}{2(2^k+1)(2^{k+1}+1)})(B(x_n-\frac{x_n}{2^k+1})-B(x_n-2\frac{x_n}{2^{k+1}+1}))\\
&&+(A(x_n-\frac{x_n}{2^k+1})-A(\frac{x_n}{2(2^k+1)(2^{k+1}+1)}))B(x_n)\\
&\le&A(\frac{x_n}{2}))o(B(x_n))+o(A(\frac{x_n}{2}))B(x_n)=o(x_n),\end{eqnarray*}
a contradiction. Then we could deduce from \eqref{08214} that \begin{eqnarray*}A(x_n-\frac{x_n}{2^{k+1}+1})-
A(x_n-2\frac{x_n}{2^{k+1}+1})=o(A(x_n)-A(x_n-\frac{x_n}{2^{k+1}+1}))=o(A(x_n))
=o(A(\frac{x_n}{2})).\end{eqnarray*} It follows that \begin{eqnarray*}A(x_n-\frac{x_n}{2^{k+1}+1})-A(x_n-\frac{x_n}{2^k+1})=o(A(\frac{x_n}{2})).
\end{eqnarray*}
Hence by \eqref{08216} we have $$A(x_n-\frac{x_n}{2^{k+1}+1})=(1+o(1))A(\frac{x_n}{2}).$$ Since $A((1-\delta)x_n)=(1+o(1))A(\frac{x_n}{2})$ for any $0<\delta <0.5$ we have $B(cx_n)=(c+o(1))B(x_n)$ for every $0<c<1$.\end{proof}\vskip2mm

We will return to the proof of Theorem \ref{thm1} (i). By Proposition 2, we know that, in any case, \begin{eqnarray*}\frac{A(y)B(y)}{y}\rightarrow1 \hskip3mm\mbox{as}\hskip3mm n\rightarrow \infty. \end{eqnarray*} Theorem \ref{thm2} (i) follows.\vskip2mm

Now we will prove Theorem \ref{thm2} (ii). We know that $A(x_n)B(x_n)=(2+o(1))x_n$ and
$$
| \{m: x_n\le m\le 2x_n, m=a+b, a\in A, b\in B, a>x_n \hskip3mm\mbox{or} \hskip3mm b>x_n\} |=o(x_n).
$$
Assume the contrary, for some $c_3x_n\le v_n\le c_4x_n$, we have $A(v_n)B(v_n)>(2+\varepsilon )x_n$ with $\varepsilon >0$ for infinitely many $n$. Then, $A(v_n)>(1+\varepsilon ')A(x_n)$ or $B(v_n)>(1+\varepsilon ')x_n$ for some $\varepsilon '>0$. Let us suppose that $A(v_n)>(1+\varepsilon ')A(x_n)$, that is $A(v_n)-A(x_n)\ge \varepsilon A(x_n)$. We proved in (i) that $B((2-c_4)x_n)=(0.5+o(1))B(x_n)$ or $B((2-c_4)x_n)=(2-c_4+o(1))B(x_n)$. It follows that $B((2-c_4)x_n)\gg B(x_n)$. Hence
\begin{eqnarray*}
o(x_n)&=&| \{m: x_n\le m\le 2x_n, m=a+b, a\in A, b\in B, a>x_n \hskip3mm\mbox{or} \hskip3mm b>x_n\} |\\
&\ge &(B(2-c_4)x_n)(A(v_n)-A(x_n))\gg A(x_n)B(x_n)=(2+o(1))x_n,
\end{eqnarray*}
a contradiction. Theorem \ref{thm2} (ii) follows.\vskip2mm

Finally we will prove Theorem \ref{thm2} (iii). By Proposition 1, with the exception of at most $\varepsilon x_n$ numbers, all numbers in $[0,2x_n]$ can be written as $a_i+b_j$ with $a_i\le x_n, b_j\le x_n$. For $i=1,2,3,4$, denote the number of elements of $A$ and $B$ in the intervals $(\frac{(i-1)x_n}{2},\frac{ix_n}{2}]$ by $A_i$ and $B_i$, respectively. Clearly,
\begin{eqnarray*}A_1B_1+A_1B_2+A_2B_1+A_2B_2\le 2x_n.\end{eqnarray*} \vskip2mm
 Noting that the sum in $A_1B_3+A_1B_4+A_2B_3+A_2B_4+A_3B_1+A_3B_2+A_3B_3+A_4B_1+A_4B_2$ belongs to $[x_n,3x_n]$, we know that
\begin{eqnarray*}A_1B_3+A_1B_4+A_2B_3+A_2B_4+A_3B_1+A_3B_2+A_3B_3+A_4B_1+A_4B_2\le x_n+\varepsilon x_n.\end{eqnarray*} \vskip2mm

Considering the difference, also by the definition of \emph{disjoint set}, with the exception of at most $\varepsilon x_n$ numbers, all numbers in $[-x_n,x_n]$ can be written as $a_i-b_j$ with $a_i\le x_n, b_j\le x_n$. Thus, \begin{eqnarray*}A_3B_4+A_4B_3+A_4B_4\le \varepsilon x_n.\end{eqnarray*}  \vskip2mm
Hence
\begin{eqnarray*}A(2x_n)B(2x_n)=(A_1+A_2+A_3+A_4)(B_1+B_3+B_3+B_4)\le 3x_n+2\varepsilon x_n.          \end{eqnarray*}
It follows that 
\begin{eqnarray*}
\limsup_{k\rightarrow \infty}\frac{A(2x_n)B(2x_n)}{2x_n}\le \frac{3}{2}.\end{eqnarray*}
\vskip2mm

On the other hand, let \begin{eqnarray}\label{08231}
&&A=\{\epsilon_0+\epsilon_2m_1m_2+\cdots+\epsilon_{2k-2}m_1\cdots m_{2k-2}, \epsilon_{2i}=0,1,\cdots,m_{2i+1}-1\},\\
&&\nonumber B=\{\epsilon_1m_1+\epsilon_3m_1m_2m_3+\cdots+\epsilon_{2k-1}m_1\cdots m_{2k-1}, \epsilon_{2i-1}=0,1,\cdots,m_{2i}-1\},
\end{eqnarray} where $m_1,m_2,\cdots$ are integers no less than two. Then $A,B$ are disjoint sets. Let \begin{eqnarray*}
y_k=(m_1-1)+(m_3-1)m_1m_2+\cdots+(m_{2k-1}-1)m_1\cdots m_{2k-2}+m_1m_2\cdots m_{2k}.\end{eqnarray*}
Then we have $A(y_n)=2m_1m_3\cdots m_{2n-1}$ and $B(y_n)=m_2m_4\cdots m_{2n}$. Let $m_{i+1}\ge im_i$ for every positive integer $i$. Then
\begin{eqnarray*}\lim_{k\rightarrow \infty}\frac{A(y_k)B(y_k)}{y_k}=2.\end{eqnarray*}
Also, $A(2y_n)=3m_1m_3\cdots m_{2n-1}$, $B(2y_n)=m_2m_4\cdots m_{2n}$
and then \begin{eqnarray*}\lim_{k\rightarrow \infty}\frac{A(2y_n)B(2y_n)}{2y_n}=\frac{3}{2}.\end{eqnarray*}
Theorem \ref{thm3} (iii) follows.
\vskip2mm

This completes the proof of Theorem \ref{thm1}. \hfill$\Box$\\

\noindent{\bf Proof of Theorem \ref{thm2}.} We will firstly prove Theorem \ref{thm2} (i). Assume the contrary, \eqref{08172} does not hold. Then there exists a positive constant $c$ such that $\frac{x_{n+1}}{x_n}\le c$ for infinitely many positive integers $n$. Since $x_{i+1}\ge 2x_i$ for all $i$, there exist infinitely many positive integers $n$ such that
\begin{eqnarray*} \frac{1}{c}x_{n+1}\le x_n\le \frac{1}{2}x_{n+1}, \hskip4mm\mbox{namely},\hskip4mm
\frac{2}{c}x_{n+1}\le 2x_n\le x_{n+1}.\end{eqnarray*} It infers from
\begin{eqnarray*}\lim_{n\rightarrow \infty}
\frac{A(x_{n+1})B(x_{n+1})}{x_{n+1}}=2\end{eqnarray*}
and Theorem \ref{thm1} (i) that \begin{eqnarray*}\lim_{n\rightarrow \infty}\frac{A(x_n)B(x_n)}{x_n}=1,\end{eqnarray*}a contradiction with \eqref{08191}.
\vskip2mm

Now we will prove Theorem \ref{thm2} (ii). We introduce the following proposition.\vskip2mm

\noindent\textbf{Proposition 3.} Define disjoint sets $A$ and $B$ as in \eqref{08231}. Write
\begin{eqnarray*}
N_k=(m_1-1)+(m_3-1)m_1m_2+\cdots+(m_{2k-1}-1)m_1\cdots m_{2k-2}.\end{eqnarray*}
Then, for $x$ with $N_k+m_1m_2\cdots m_{2k}\le x\le 2m_1m_2\cdots m_{2k-1}+m_1m_2\cdots m_{2k}$, we have
\begin{eqnarray*}\label{08181} \frac{A(x)B(x)}{x}\ge \frac{2}{1+\frac{2}{m_{2k}}}.
\end{eqnarray*}
\vskip2mm

We now return to the proof of Theorem \ref{thm2} (ii). Suppose we have already found $m_1,m_2,\cdots,m_{2k-2}$ such that
\begin{eqnarray*}N_{k-1}+m_1m_2\cdots m_{2k-2}\le x_i\le 2m_1m_2\cdots m_{2k-3}+m_1m_2\cdots m_{2k-2}.\end{eqnarray*} Now, for $x_{i+1}$, we need to find $m_{2k-1}$ and $m_{2k}$ such that
\begin{eqnarray*}\label{08182}N_{k}+m_1m_2\cdots m_{2k}\le x_{i+1}\le 2m_1m_2\cdots m_{2k-1}+m_1m_2\cdots m_{2k}.\end{eqnarray*} Let $m_{2k-1}=2$. Since the length of the interval
$[\frac{x_{i+1}-2m_1m_2\cdots m_{2k-1}}{m_1m_2\cdots m_{2k-1}},
\frac{x_{i+1}-N_k}{m_1m_2\cdots m_{2k-1}}]$ is larger than 1, there exists an integer (write it as $m_{2k}$) such that
\begin{eqnarray}\label{08184}
\frac{x_{i+1}-2m_1m_2\cdots m_{2k-1}}{m_1m_2\cdots m_{2k-1}}\le m_{2k}
\le \frac{x_{i+1}-N_k}{m_1m_2\cdots m_{2k-1}}.\end{eqnarray}
Namely,
\begin{eqnarray*}
N_k+m_1m_2\cdots m_{2k}\le x_{i+1}\le 2m_1m_2\cdots m_{2k-1}+m_1m_2\cdots m_{2k}.\end{eqnarray*}
It infers from Proposition 3 that
\begin{eqnarray*}\frac{A(x_{i+1})B(x_{i+1})}{x_{i+1}}
\ge \frac{2}{1+\frac{2}{m_{2k}}}.\end{eqnarray*}
By \eqref{08172} and \eqref{08184} we know that $m_{2k}\rightarrow\infty$. Thus, \begin{eqnarray*}
\lim_{n\rightarrow \infty}\frac{A(x_n)B(x_n)}{x_n}=2.\end{eqnarray*}
\vskip2mm

This completes the proof of Theorem \ref{thm2}. \hfill$\Box$\\

The main idea in the proof of Theorem \ref{thm3} is from \cite[Theorem 3.3]{Erdos2}. We give the details for the sake of completeness. \vskip2mm

\noindent{\bf Proof of Theorem \ref{thm3}.} If $IN=0$, then \eqref{thm201} follows from the fact that $SP\le 2$. Now we only need to consider the case $IN>0$. For any $\eta$ with $0<\eta<IN$, we have $A(x)\ge(IN-\eta)\sqrt{x}$ and $B(x)\ge(IN-\eta)\sqrt{x}$ for large enough $x$. It follows from $A(x)B(x)\le 2x$ that
\begin{eqnarray}\label{03191}
B(x)\le\frac{2}{IN-\eta}\sqrt{x} \hskip3mm\mbox{and}\hskip3mm A(x)\le\frac{2}{IN-\eta}\sqrt{x} \hskip3mm\mbox{for large} \hskip3mm x.\end{eqnarray}\vskip2mm

Let $\varepsilon$ be small enough. Take an $x$ such that
\begin{eqnarray*}A(2x)B(2x)\ge(2SP-\varepsilon)x.\end{eqnarray*}
Thus, with the exception of at most $(4-2SP+\varepsilon)x$ numbers, all numbers in $[0,4x]$ can be written as $a_i+b_j$ with $a_i\le 2x, b_j\le 2x$. Let $$A_1=|A\cap [0,x]|, A_2=|A\cap (x,2x]|, B_1=|B\cap [0,x]|, B_2=|B\cap (x,2x]|.$$ Hence
\begin{eqnarray}\label{08233}
A_1B_1+A_2B_2\ge(2SP-2-\varepsilon)x\end{eqnarray}
and also \begin{eqnarray*}\label{03194}
A_1B_1\ge(2SP-3-\varepsilon)x,\hskip3mm A_2B_2\ge(2SP-3-\varepsilon)x.\end{eqnarray*}
It follows from $A_1B_1+A_2B_2\le 2x$ that \begin{eqnarray*}\label{08232}
A_1B_1\le(5-2SP+\varepsilon)x,\hskip3mm A_2B_2\le(5-2SP+\varepsilon)x.\end{eqnarray*}\vskip2mm

Take $$d=\frac{(IN-\eta)^4}{16}.$$ Then $0<d<1$. Let
$$A'=|A\cap [dx,x]|, A^{*}=|A\cap (x,(1+d)x]|, B'=|B\cap [dx,x]|, B^{*}=|B\cap (x,(1+d)x]|.$$
Considering the difference, by \eqref{08233} we know that at most $(4-2SP+\varepsilon)x$ other pairs satisfy $|a_i-b_j|\le x$. Thus, \begin{eqnarray}\label{03193}
A'B^{*}+A^{*}B'\le (4-2SP+\varepsilon)x.\end{eqnarray}
By \eqref{03191} we have
$$A'=A(x)-A(dx)\ge(IN-\eta)\sqrt{x}-\frac{2}{IN-\eta}\sqrt{dx}=\frac{IN-\eta}{2}\sqrt{x},$$(similarly, $B'>\frac{IN-\eta}{2}\sqrt{x}$.) It infers from \eqref{03193} that
$$A^{*}+B^{*}\le \frac{2(4-2SP+\varepsilon)}{IN-\eta}\sqrt{x}.$$ Write
$M=\max\{A^{*},B^{*}\}$. Then
$$M\le \frac{2(4-2SP+\varepsilon)}{IN-\eta}\sqrt{x}.$$
By \eqref{08232}, $$A((1+d)x)B((1+d)x)\ge(1+d)x-(4-2SP+\varepsilon)x=(2SP+d-3-\varepsilon)x$$
and $$A((1+d)x)B((1+d)x)=(A_1+A^{*})(B_1+B^{*}),$$ we have
$$A_1B^{*}+A^{*}B_1+A^{*}B^{*}>(4SP+d-8-2\varepsilon)x.$$
On the other hand,
$$A_1B^{*}+A^{*}B_1+A^{*}B^{*}\le M(A_1+B_1)+M^2
<\frac{8\sqrt{x}}{IN-\eta}M\le \frac{8}{IN-\eta}\frac{2(4-2SP+\varepsilon)}{IN-\eta}x.$$
Thus,
$$4SP+d-8-2\varepsilon\le\frac{16(4-2SP+\varepsilon)}{(IN-\eta)^2}.$$
As $\eta\rightarrow0$ and $\varepsilon\rightarrow0$, we have $$\frac{IN^6}{IN^2+8}\le 64(2-SP).$$ \eqref{thm201} follows.\vskip2mm

This completes the proof of Theorem \ref{thm3}. \hfill$\Box$

\section{Acknowledgments}
The first author sincerely thanks Prof. Yong-Gao Chen for his kind help and concern all the time.

\end{document}